# Numerical Solutions of Fredholm Integral Equations Using Bernstein Polynomials

**A. Shirin[1], M. S. Islam[2]**

[1]Institute of Natural Sciences, United International University, Dhaka-1209, Bangladesh

[2]Department of Mathematics, University of Dhaka, Dhaka-1000, Bangladesh

### Abstract

In this paper, Bernstein piecewise polynomials are used to solve the integral equations numerically. A matrix formulation is given for a non-singular linear Fredholm Integral Equation by the technique of Galerkin method. In the Galerkin method, the Bernstein polynomials are exploited as the linear combination in the approximations as basis functions. Examples are considered to verify the effectiveness of the proposed derivations, and the numerical solutions guarantee the desired accuracy.

## 1. Introduction

In the survey of solutions of integral equations, a large number of analytical but a few approximate methods for solving numerically various classes of integral equations [1, 2] are available. Since the piecewise polynomials are differentiable and integrable, the Bernstein polynomials [4, 5] are defined on an interval to form a complete basis over the finite interval. Moreover, these polynomials are positive and their sum is unity. For these advantages, Bernstein polynomials have been used to solve second order linear and nonlinear differential equations, which are available in the literature, e.g. Bhatti and Bracken [7]. Very recently, Mandal and Bhattacharya [6] have attempted to solve integral equations numerically using Bernstein polynomials, but they obtained the results in terms of finite series solutions. In contrast to this, we solve the linear Fredholm integral equation by exploiting very well known Galerkin method [3] and Bernstein polynomials are used as trial functions in the basis. For this, we give a short introduction of Bernstein polynomials first. Then we derive a matrix formulation by the technique of Galerkin method. To verify our formulation we consider three examples, in which we obtain exact solutions for two examples even using a few and lower order polynomials. On the other hand, the last example shows an excellent agreement of accuracy compared to exact solution, which confirms the convergence. All the computations are performed using *MATHEMATICA*.

---

[2] *Corresponding author*:mdshafiqul@yahoo.com

## 2. Bernstein Polynomials

The general form of the Bernstein polynomials [4-7] of $n$th degree over the interval $[a,b]$ is defined by

$$B_{i,n}(x) = \binom{n}{i} \frac{(x-a)^i (b-x)^{n-i}}{(b-a)^n}, \quad a \le x \le b, \quad i = 0,1,2,\cdots,n \quad (1)$$

Note that each of these $n+1$ polynomials having degree $n$ satisfies the following properties:

i) $B_{i,n}(x) = 0$, if $i < 0$ or $i > n$, ii) $\sum_{i=0}^{n} B_{i,n}(x) = 1$;

iii) $B_{i,n}(a) = B_{i,n}(b) = 0$, $1 \le i \le n-1$

Using *MATHEMATICA* code, the first 11 Bernstein polynomials of degree ten over the interval $[a, b]$, are given below:

$B_{0,10}(x) = (b-x)^{10}/(b-a)^{10}$

$B_{1,10}(x) = 10(b-x)^9(x-a)/(b-a)^{10}$

$B_{2,10}(x) = 45(b-x)^8(x-a)^2/(b-a)^{10}$

$B_{3,10}(x) = 120(b-x)^7(x-a)^3/(b-a)^{10}$

$B_{4,10}(x) = 210(b-x)^6(x-a)^4/(b-a)^{10}$

$B_{5,10}(x) = 252(b-x)^5(x-a)^5/(b-a)^{10}$

$B_{6,10}(x) = 210(b-x)^4(x-a)^6/(b-a)^{10}$

$B_{7,10}(x) = 120(b-x)^3(x-a)^7/(b-a)^{10}$

$B_{8,10}(x) = 45(b-x)^2(x-a)^8/(b-a)^{10}$

$B_{9,10}(x) = 10(b-x)(x-a)^9/(b-a)^{10}$

$B_{10,10}(x) = (x-a)^{10}/(b-a)^{10}$

Now the first six polynomials over $[0, 1]$ are shown in Fig. 1($a$), and the remaining five polynomials are shown in Fig. 1($b$).

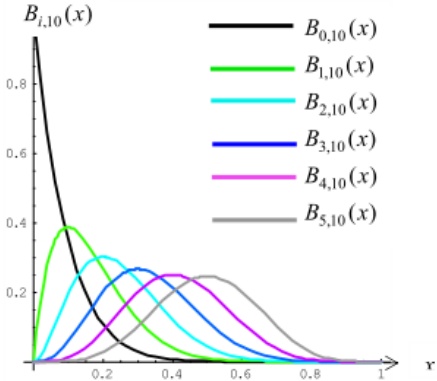 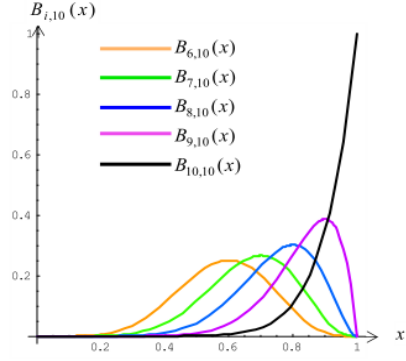

Fig. 1(a): Graph of first 6 Bernstein polynomials over [0, 1]  Fig. 1(b): Graph of last 5 Bernstein polynomials over [0, 1]

## 3. Formulation of Integral Equation in Matrix Form

Consider a general linear Fredholm integral equation (FIE) of second kind [1, 2] is given by

$$a(x)\phi(x) + \lambda \int_a^b k(t, x)\phi(t)dt = f(x) , \; a \leq x \leq b \qquad (2)$$

where $a(x)$ and $f(x)$ are given functions, $k(t, x)$ is the kernel, and $\phi(x)$ is the unknown function or exact solution of (2), which is to be determined.

Now we use the technique of Galerkin method [Lewis, 3] to find an approximate solution $\tilde{\phi}(x)$ of (2). For this, we assume that

$$\tilde{\phi}(x) = \sum_{i=0}^{n} a_i \, B_{i,n}(x) \qquad (3)$$

where $B_{i,n}(x)$ are Bernstein polynomials (basis) of degree $i$ defined in eqn. (1), and $a_i$ are unknown parameters, to be determined. Substituting (3) into (2), we obtain

$$a(x)\sum_{i=0}^{n} a_i \, B_{i,n}(x) + \lambda \int_a^b \left[ k(t,x) \sum_{i=0}^{n} a_i \, B_{i,n}(t) \right] dt = f(x)$$

or,

$$\sum_{i=0}^{n} a_i \left[ a(x) \, B_{i,n}(x) + \lambda \int_a^b k(t,x) \, B_{i,n}(t) dt \right] = f(x) \qquad (4)$$

Then the Galerkin equations [Lewis, 3] are obtained by multiplying both sides of (3) by $B_{j,n}(x)$ and then integrating with respect to $x$ from $a$ to $b$, we have

$$\sum_{i=0}^{n} a_i \left[ \int_a^b \left[ a(x) B_{i,n}(x) + \lambda \int_a^b k(t,x) B_{i,n}(t) dt \right] B_{j,n}(x) dx \right] = \int_a^b B_{j,n}(x) f(x) dx, \quad j=0,1,\cdots,n$$

Since in each equation, there are three integrals. The inner integrand of the left side is a function of $x$ and $t$, and is integrated with respect to $t$ from $a$ to $b$. As a result the outer integrand becomes a function of $x$ only and integration with respect to $x$ yields a constant. Thus for each $j$ $(=0,1,\ldots,n)$ we have a linear equation with $n+1$ unknowns $a_i$ ($i=0,1,\ldots,n$). Finally (5a) represents the system of $n+1$ linear equations in $n+1$ unknowns, are given by

$$\sum_{i=0}^{n} a_i C_{i,j} = F_j, \quad j=0,1,2,\cdots,n \tag{5a}$$

where

$$C_{i,j} = \int_a^b \left[ a(x) B_{i,n}(x) + \lambda \int_a^b k(t,x) B_{i,n}(t) dt \right] B_{j,n}(x) dx, \qquad i,j=0,1,2,\ldots,n. \tag{5b}$$

$$F_j = \int_a^b B_{i,n}(x) f(x) dx, \quad j=0,1,2,\cdots,n \tag{5c}$$

Now the unknown parameters $a_i$ are determined by solving the system of equations (5), and substituting these values of parameters in (3), we get the approximate solution $\tilde{\phi}(x)$ of the integral equation (2). The absolute error $E$ for this formulation is defined by

$$E = \left| \frac{\varphi(x) - \tilde{\varphi}(x)}{\varphi(x)} \right|.$$

## 4. Numerical Examples

In this section, we explain three integral equations which are available in the existing literatures [1, 2, 6]. For each example we find the approximate solutions using different number of Bernstein polynomials.

**Example 1:** We consider the FIE of 2nd kind given by [6]

$$\phi(x) - \int_{-1}^{1} (xt + x^2 t^2) \phi(t) dt = 1, \quad -1 \le x \le 1 \tag{6}$$

having the exact solution, $\phi(x) = 1 + \dfrac{10}{9} x^2$

Using the formulation described in the previous section, the equations (5) lead us, respectively,

$$\sum_{i=0}^{n} a_i\, C_{i,j} = F_j, \qquad j = 0,1,2,\ldots,n, \tag{7a}$$

$$C_{i,j} = \int_{-1}^{1} B_{i,n}(x) B_{j,n}(x)\,dx - \int_{-1}^{1}\left[\int_{-1}^{1}(xt + x^2 t^2) B_{i,n}(t)\,dt\right] B_{j,n}(x)\,dx, \quad i,j = 0,1,\ldots,n \tag{7b}$$

$$F_j = \int_{-1}^{1} B_{j,n}(x)\,dx, \qquad j = 0,1,2,\ldots,n \tag{7c}$$

Solving the system (7) for $n = 3$, the values of the parameters are:

$$a_0 = \dfrac{19}{9},\ a_1 = \dfrac{17}{27},\ a_2 = \dfrac{17}{27},\ a_3 = \dfrac{19}{9}$$

Substituting into (3) and for $n \geq 3$, the approximate solution is,

$$\tilde{\phi}(x) = 1 + \dfrac{10}{9} x^2$$

which is the exact solution.

**Example 2:** Now we consider another FIE of 2nd kind given by [6]

$$\phi(x) - \int_{-1}^{1}(x^4 - t^4)\phi(t)\,dt = x, \quad -1 \leq x \leq 1 \tag{8}$$

having the exact solution $\phi(x) = x$

Proceeding as the example 1, the system of equations becomes as

$$\sum_{i=0}^{n} a_i C_{i,j} = F_j,\quad j = 0,1,\ldots,n \tag{9a}$$

where,

$$C_{i,j} = \int_{-1}^{1} B_{i,n}(x) B_{i,n}(x)\,dx - \int_{-1}^{1}\left[\int_{-1}^{1}(x^4 - t^4) B_{i,n}(t)\,dt\right] B_{j,n}(x)\,dx,\quad i,j = 0,1,\ldots,n \tag{9b}$$

$$F_j = \int_{-1}^{1} x B_{j,n}(x)dx, \qquad j = 0,1,\ldots,n \qquad (9c)$$

Now solving the system (9) for $n = 3$, the values of the parameters, $a_i$ are:

$$a_0 = -1, \quad a_1 = -\frac{1}{3}, \quad a_2 = \frac{1}{3}, \quad a_3 = 1.$$

and the approximate solution, for $n \geq 3$, is $\tilde{\phi}(x) = x$ which is the exact solution.

**Example 3:** Consider another FIE of 2nd kind given by [1, pp 213]

$$\phi(x) - \int_0^1 (tx^2 + xt^2)\phi(t)dt = x, \qquad 0 \leq x \leq 1 \qquad (10)$$

having the exact solution $\phi(x) = \frac{180}{119}x + \frac{80}{119}x^2$

Proceeding as the previous examples, the system of equations becomes:

$$C_{i,j} = \int_0^1 B_{i,n}(x) B_{j,n}(x)dx - \int_0^1 \left[ \int_0^1 (tx^2 + xt^2) B_{i,n}(t)dt \right] B_{j,n}(x)dx, \quad i,j = 0,1,\ldots,n \quad (11a)$$

$$F_j = \int_0^1 x B_{j,n}(x)dx, \qquad j = 0,1,\ldots,n \qquad (11b)$$

For $n = 3$, solving system (11), the values of the parameters ($a_i$) are:

$$a_0 = 0, \quad a_1 = \frac{260}{119}, \quad a_2 = \frac{80}{119}, \quad a_3 = 0$$

and the approximate solution is $\tilde{\phi}(x) = \frac{180}{119}x + \frac{80}{119}x^2$, which is the exact solution.

**Example 4:** Consider another FIE of 2nd kind given by [2, pp 124]

$$\phi(x) - \int_0^1 2e^x e^t \phi(t)dt = e^x, \qquad 0 \leq x \leq 1, \qquad (12)$$

having the exact solution $\phi(x) = \frac{e^x}{2 - e^2}$.

Since the equations (5b) and (5c) are of the form

$$C_{i,j} = \int_0^1 B_{i,n}(x) B_{j,n}(x) dx - \int_0^1 \left[ \int_0^1 2e^x e^t B_{i,n}(t) dt \right] B_{j,n}(x) dx \quad i, j = 0,1,2,\ldots,n \quad (13a)$$

$$F_j = \int_0^1 e^x B_{j,n}(x) dx, \quad j = 0,1,\ldots,n \quad (13b)$$

For $n = 3,4,5$, and 6, the approximate solutions are, respectively

$$\tilde{\phi}(x) = -0.185387 - 0.188957\ x - 0.078167\ x^2 - 0.051702\ x^3$$

$$\tilde{\phi}(x) = -0.185571 - 0.185273\ x - 0.0947442\ x^2 - 0.0259153\ x^3 - 0.0128933\ x^4$$

$$\tilde{\phi}(x) = -0.185561 - 0.18558\ x - 0.0926328\ x^2 - 0.0315503\ x^3 - 0.00655026\ x - 0.00253837\ x^5$$

$$\tilde{\phi}(x) = -0.185561 - 0.18556\ x - 0.0925932\ x^2 - 0.0308581 x^3 - 0.00791665\ x^4$$
$$- 0.0012903\ x^5 - 0.000427923 x^6$$

Plot of absolute difference, the error *E* between exact and approximate solutions, is depicted in Fig. 2 for various values of *n*. Observe that the minimum order of accuracies are $10^{-4}$, $10^{-6}$, $10^{-7}$, and $10^{-8}$, respectively, with 4, 5, 6 and 7 Bernstein polynomials. This confirms us that we if increase the number of polynomials, the accuracy also increases. Now the approximate solutions, exact solutions, and the error *E*, between exact and the approximate solutions at various points of the domain are displayed in Table 1.

Table 1. Numerical solutions at various points and corresponding absolute errors of the example 4.

| x | Exact Solutions | Approximate Solutions | Absolute Relative Error, $E$ | Approximate Solutions | Absolute Relative Error, $E$ |
|---|---|---|---|---|---|
| | | Polynomials used 4 | | Polynomials used 5 | |
| 0.0 | -0.1855612526 | -0.1853868426 | 0.000940 | -0.1855710276 | 0.0000526782 |
| 0.1 | -0.2050768999 | -0.2051159200 | 0.000190 | -0.2050729953 | 0.0000190395 |
| 0.2 | -0.2266450257 | -0.2267185494 | 0.000324 | -0.2266433896 | $7.218900 \times 10^{-6}$ |
| 0.3 | -0.2504814912 | -0.2505049431 | 0.000094 | -0.2504841183 | 0.0000104884 |
| 0.4 | -0.2768248595 | -0.2767853131 | 0.000143 | -0.2768280333 | 0.0000114649 |
| 0.5 | -0.3059387842 | -0.3058698717 | 0.000225 | -0.3059389305 | $4.784034 \times 10^{-7}$ |
| 0.6 | -0.3381146470 | -0.3380688310 | 0.000136 | -0.3381115499 | $9.159888 \times 10^{-6}$ |
| 0.7 | -0.3736744748 | -0.3736924032 | 0.000048 | -0.3736715753 | $7.759438 \times 10^{-6}$ |
| 0.8 | -0.4129741624 | -0.4130508005 | 0.000186 | -0.4129756348 | $3.565315 \times 10^{-6}$ |
| 0.9 | -0.4564070342 | -0.4564542350 | 0.000103 | -0.4564113003 | $9.347079 \times 10^{-6}$ |
| 1.0 | -0.5044077810 | -0.5042129189 | 0.000386 | -0.5043970878 | 0.0000211995 |
| | | Polynomials used 6 | | Polynomials used 7 | |
| 0.0 | -0.1855612526 | -0.1855612526 | $1.358183 \times 10^{-6}$ | -0.1855612694 | 0.0000415721 |
| 0.1 | -0.2050768999 | -0.2050768100 | $5.311130 \times 10^{-7}$ | -0.2050768958 | 0.0000127221 |
| 0.2 | -0.2266450257 | -0.2266450257 | $5.268384 \times 10^{-7}$ | -0.2266450312 | $3.763352 \times 10^{-6}$ |
| 0.3 | -0.2504814912 | -0.2504814912 | $4.307390 \times 10^{-7}$ | -0.2504814909 | $6.065680 \times 10^{-6}$ |
| 0.4 | -0.2768248595 | -0.2768248596 | $2.997736 \times 10^{-7}$ | -0.2768248544 | $3.529983 \times 10^{-6}$ |
| 0.5 | -0.3059387842 | -0.3059387842 | $5.281912 \times 10^{-7}$ | -0.3059387842 | $7.651093 \times 10^{-6}$ |
| 0.6 | -0.3381146470 | -0.3381146470 | $3.685661 \times 10^{-8}$ | -0.3381146522 | $1.444431 \times 10^{-6}$ |
| 0.7 | -0.3736744748 | -0.3736744748 | $4.654666 \times 10^{-7}$ | -0.3736744750 | $5.945203 \times 10^{-6}$ |
| 0.8 | -0.4129741624 | -0.4129741624 | $1.890301 \times 10^{-7}$ | -0.4129741564 | $3.339474 \times 10^{-6}$ |
| 0.9 | -0.4564070342 | -0.4564070342 | $5.100054 \times 10^{-7}$ | -0.4564070387 | $5.480605 \times 10^{-6}$ |
| 1.0 | -0.5044077810 | -0.5044077810 | $1.161812 \times 10^{-6}$ | -0.5044077618 | 0.0000149584 |

Table 1. Numerical solutions at various points and corresponding absolute errors of the example 4.

| x | Exact Solutions | Approximate Solutions | Absolute Relative Error, $E$ | Approximate Solutions | Absolute Relative Error, $E$ |
|---|---|---|---|---|---|
| | | Polynomials used 4 | | Polynomials used 5 | |
| 0.0 | -0.1855612526 | -0.1853868426 | 0.000940 | -0.1855710276 | 0.0000526782 |
| 0.1 | -0.2050768999 | -0.2051159200 | 0.000190 | -0.2050729953 | 0.0000190395 |
| 0.2 | -0.2266450257 | -0.2267185494 | 0.000324 | -0.2266433896 | $7.218900 \times 10^{-6}$ |
| 0.3 | -0.2504814912 | -0.2505049431 | 0.000094 | -0.2504841183 | 0.0000104884 |
| 0.4 | -0.2768248595 | -0.2767853131 | 0.000143 | -0.2768280333 | 0.0000114649 |
| 0.5 | -0.3059387842 | -0.3058698717 | 0.000225 | -0.3059389305 | $4.784034 \times 10^{-7}$ |
| 0.6 | -0.3381146470 | -0.3380688310 | 0.000136 | -0.3381115499 | $9.159888 \times 10^{-6}$ |
| 0.7 | -0.3736744748 | -0.3736924032 | 0.000048 | -0.3736715753 | $7.759438 \times 10^{-6}$ |
| 0.8 | -0.4129741624 | -0.4130508005 | 0.000186 | -0.4129756348 | $3.565315 \times 10^{-6}$ |
| 0.9 | -0.4564070342 | -0.4564542350 | 0.000103 | -0.4564113003 | $9.347079 \times 10^{-6}$ |
| 1.0 | -0.5044077810 | -0.5042129189 | 0.000386 | -0.5043970878 | 0.0000211995 |
| | | Polynomials used 6 | | Polynomials used 7 | |
| 0.0 | -0.1855612526 | -0.1855612526 | $1.358183 \times 10^{-6}$ | -0.1855612694 | 0.0000415721 |
| 0.1 | -0.2050768999 | -0.2050768100 | $5.311130 \times 10^{-7}$ | -0.2050768958 | 0.0000127221 |
| 0.2 | -0.2266450257 | -0.2266450257 | $5.268384 \times 10^{-7}$ | -0.2266450312 | $3.763352 \times 10^{-6}$ |
| 0.3 | -0.2504814912 | -0.2504814912 | $4.307390 \times 10^{-7}$ | -0.2504814909 | $6.065680 \times 10^{-6}$ |
| 0.4 | -0.2768248595 | -0.2768248596 | $2.997736 \times 10^{-7}$ | -0.2768248544 | $3.529983 \times 10^{-6}$ |
| 0.5 | -0.3059387842 | -0.3059387842 | $5.281912 \times 10^{-7}$ | -0.3059387842 | $7.651093 \times 10^{-6}$ |
| 0.6 | -0.3381146470 | -0.3381146470 | $3.685661 \times 10^{-8}$ | -0.3381146522 | $1.444431 \times 10^{-6}$ |
| 0.7 | -0.3736744748 | -0.3736744748 | $4.654666 \times 10^{-7}$ | -0.3736744750 | $5.945203 \times 10^{-6}$ |
| 0.8 | -0.4129741624 | -0.4129741624 | $1.890301 \times 10^{-7}$ | -0.4129741564 | $3.339474 \times 10^{-6}$ |
| 0.9 | -0.4564070342 | -0.4564070342 | $5.100054 \times 10^{-7}$ | -0.4564070387 | $5.480605 \times 10^{-6}$ |
| 1.0 | -0.5044077810 | -0.5044077810 | $1.161812 \times 10^{-6}$ | -0.5044077618 | 0.0000149584 |

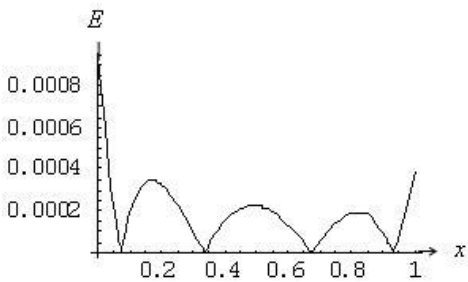

Fig. 2a. Error $E$, using 4 polynomials

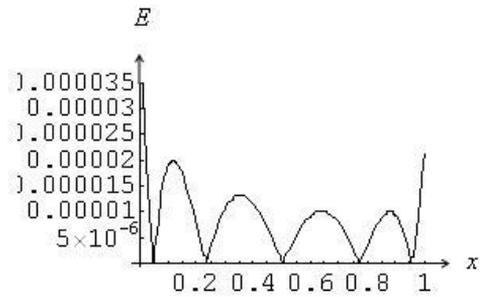

Fig. 2b. Error $E$, using 5 polynomials

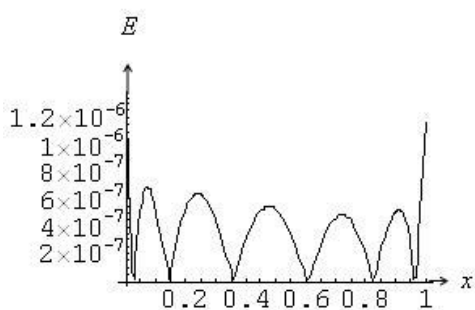 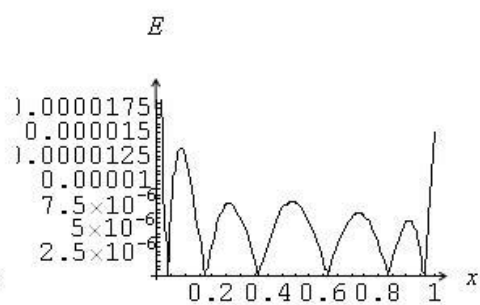

Fig. 2c. Error E, using 6 polynomials    Fig. 2d. Error E, using 7 polynomials

## 5. Conclusion

We have considered the integral equations to solve numerically. We have obtained the approximate solution of the unknown function by the well known Galerkin method using Bernstein polynomials as trial functions. We have verified the derived formulas with the appropriate numerical examples. In this context we may note that the numerical solutions coincide with the exact solutions even a few of the polynomials are used in the approximation.

**References**       **Ref format/style not OK**